\documentclass[conf]{new-aiaa}
\usepackage[utf8]{inputenc}
\usepackage[T1]{fontenc}
\usepackage[table]{xcolor}
\usepackage{amsmath} 
\usepackage{graphicx,subcaption}
\usepackage{boldline}
\usepackage{tikz}
\usepackage{authblk}
\usepackage{secdot}
\usepackage{hyperref}
\hypersetup{
    colorlinks=true,
    linkcolor=blue,
    citecolor=blue,
    filecolor=magenta,      
    urlcolor=cyan,
}
\usepackage{cleveref}
\usepackage{todonotes}
\usepackage[linesnumbered,ruled,vlined]{algorithm2e}

\title{Tensor-Train Operator Inference}
\author{Engin Danis\footnote{Assistant Professor, Department of Mechanical and Aerospace Engineering, AIAA Member}\footnote{Theoretical Division at Los Alamos National Labarotory}}
\affil{University of Missouri, Columbia, MO 65211}
\author{and \\Duc Truong\footnote{Scientist, Theoretical Division at Los Alamos National Labarotory} and Kim {\O.} Rasmussen\footnote{Scientist, Theoretical Division at Los Alamos National Labarotory} and Boian S. Alexandrov\footnote{Senior Scientist, Theoretical Division at Los Alamos National Labarotory} }
\affil{Los Alamos National Laboratory, Los Alamos, NM 87545}
\begin{document}

\maketitle

\begin{abstract}
In this study, we present a tensor-train framework for nonintrusive operator inference aimed at learning discrete operators and using them to predict solutions of physical governing equations. Our framework comprises three approaches: full-order tensor-train operator inference, full-order quantized tensor-train operator inference, and reduced-order tensor-train operator inference. In each case, snapshot data is represented in tensor-train format—either through compression or cross interpolation—enabling the efficient handling of extremely large datasets with significantly reduced computational effort compared to standard methods. The effectiveness of each approach is demonstrated through numerical experiments related to Computational Fluid Dynamics and benchmarked against the standard reduced-order operator inference method, highlighting the advantages of the tensor-train representations in both accuracy and scalability.
\end{abstract}

% \tableofcontents

% {\color{red} TO DO:
% \begin{itemize}
%     \item introduction, lit review, etc.
%     \item go over the methods section and improve what's already written for the extended abstract.
%     \item Add new sections for TT implementation into the methods section
%     \item discuss our cross interpolation approach that avoids loading the complete snapshot to the memory.
%     \item discuss how we implement regularization in the full-order TT/QTT OpInf: recall the standard OpInf can use two different reg factors for the convection and diff operators but full-order TT can only use one reg factor for both (see the pseudo inverse implementation in the code) 
%     \item discuss how Full-order TT/QTT uses tt\_tensor for initial conversion while TTROM uses cross interpolation.
% \end{itemize}}
\section{Introduction}
Reduced-order modeling provides an efficient framework for accelerating numerical simulations of complex physical systems by projecting dynamics onto low-dimensional subspaces \cite{benner2015survey,quarteroni2015reduced,brunton2022data}. Specifically, Operator Inference (OpInf) has emerged as a data-driven approach within this context \cite{peherstorfer2016opinf,swischuk2019projection}. OpInf enables direct construction of reduced operators from snapshot data without explicit access to the governing equations. Assuming a polynomial structure in the dynamics, OpInf works through a least-squares problem to infer the operators that approximate system dynamics. 

While effective for moderate-scale problems, the standard formulation of OpInf encounters computational bottlenecks for high-dimensional problems. The approach typically relies on Proper Orthogonal Decomposition (POD), which involves a singular value decomposition (SVD) of potentially massive snapshot matrices \cite{li2022enhanced}. For large-scale systems, the associated memory and computational costs become prohibitive. Moreover, the inclusion of nonlinear terms increases the complexity, as their reduced representations scale polynomially with the number of retained modes, increasing the overall expense. 

These limitations are a manifestation of the curse of dimensionality in scientific computing \cite{bellman1966dynamic} and motivates the use of tensor-based methods for data compression and more efficient computing. The Tensor-Train (TT) decomposition offers an attractive format for representing high-dimensional tensors in a low-rank form. By expressing a high-dimensional array as a sequence of low-dimensional cores connected via matrix products, the TT format significantly reduce storage requirements. Importantly, many algebraic operations can be performed directly within the TT representation avoiding to expand the tensors in full. 

Here we propose a tensorized extension of OpInf, termed TT-OpInf, designed to overcome the aforementioned computational barriers. Two variants are developed for different use cases. A full-order TT-OpInf formulation recasts the operator learning process entirely in tensor format. For structured grids, this tensorization naturally reflects the spatial arrangement of data; for unstructured problems, the Quantized Tensor-Train (QTT) approach \cite{khoromskij2011d}—based on prime factorization—yields a compatible tensor structure suitable for TT operations. Finally, a reduced-order variant of TT-OpInf is introduced as an alternative to POD-based dimensionality reduction. Rather than computing an SVD of the full snapshot matrix, the snapshot tensor is converted directly into TT format using TT-SVD \cite{oseledets2010tt} or TT cross-approximation algorithms \cite{savostyanov2011fast}.The latter approach is particularly noteworthy as it constructs the TT decomposition using selected entries only, avoiding the need to store the entire dataset. 

Central to our method is the development of least-squares solvers that operate fully within the TT framework. We present algorithms for computing regularized pseudoinverses of TT tensors and for performing necessary contractions to extract reduced operator representations. All computations are maintained in compressed format, ensuring scalability with system size.
To illustrate the capabilities of the method, we conduct numerical experiments on three benchmark problems of increasing complexity: the two-dimensional heat equation, the viscous Burgers' equation, and compressible laminar flow over a cylinder. These examples demonstrate the accuracy, memory efficiency, and computational advantages of TT-OpInf, highlighting the trade-offs between full-order and reduced-order formulations.

The remainder of this paper is organized as follows. Section II reviews the standard OpInf framework and the fundamentals of the TT format. Section III introduces the TT-OpInf methodology in both full-order and reduced-order settings. Section IV presents numerical results, and Section V concludes with a summary and discussion of future directions.

\section{Methods}
In this section, we first review the the standard OpInf method and Tensor-Train method. Thereafter the tensorization of the full-order OpInf method is presented along with the TT implementation. Finally, we present the reduced order TT-OpInf method.

\subsection{Operator Inference (OpInf)}\label{sec:opinf-standard}

Operator inference (OpInf) is a data-driven modeling technique for learning reduced-order models (ROMs) from time-resolved simulation or experimental data. It is particularly effective for systems governed by ODEs with polynomial structure \cite{peherstorfer2016opinf, kramer2024review}.

\subsubsection{Problem Formulation} \label{sec:opinf-formulation}

We begin with a set of full-state snapshots \( \{ x_k \in \mathbb{R}^N \}_{k=1}^K \), sampled at discrete time instances \( \{t_k\}_{k=1}^K \). These are assembled into the snapshot matrix
\[
X = [x_1 \;\; x_2 \;\; \cdots \;\; x_K]^T \in \mathbb{R}^{K \times N},
\]
with corresponding time derivatives
\[
\dot{X} = [\dot{x}_1 \;\; \dot{x}_2 \;\; \cdots \;\; \dot{x}_K]^T \in \mathbb{R}^{K \times N}.
\]
We assume that the system follows a quadratic dynamical model:
\begin{equation}
\label{eq:full-ode}
\frac{dx}{dt} = A x + F(x \otimes x) + B u(t),
\end{equation}
where:
\begin{itemize}
    \item \( A \in \mathbb{R}^{N \times N} \): linear operator,
    \item \( F \in \mathbb{R}^{N \times N^2} \): quadratic operator (with Kronecker structure),
    \item \( B \in \mathbb{R}^{N \times p} \): input operator,
    \item \( u(t) \in \mathbb{R}^p \): known time-dependent input.
\end{itemize}
To reduce dimensionality, we perform proper orthogonal decomposition (POD) on \( X \) and retain the first \( n \) right singular vectors to form the basis matrix \( V_n \in \mathbb{R}^{N \times n} \). The reduced coordinates are:
\[
\hat{X} = X V_n \in \mathbb{R}^{K \times n}, \quad \hat{\dot{X}} = \dot{X} V_n \in \mathbb{R}^{K \times n}.
\]

We assume the reduced-order model takes the form:
\begin{equation}
\label{eq:reduced-ode}
\frac{d\hat{x}}{dt} = \hat{A} \hat{x} + \hat{F}(\hat{x} \otimes \hat{x}) + \hat{B} u(t),
\end{equation}
where \( \hat{A} \in \mathbb{R}^{n \times n} \), \( \hat{F} \in \mathbb{R}^{n \times n^2} \), and \( \hat{B} \in \mathbb{R}^{n \times p} \) are the reduced operators to be learned.

To learn these operators from data, we define the design matrix \( D \in \mathbb{R}^{K \times (n + p + n^2)} \) and the target matrix \( R \in \mathbb{R}^{K \times n} \) as:
\begin{align}
D &= \begin{bmatrix} \hat{X} & U & \hat{X}^{(2)} \end{bmatrix}, \\
\hat{X}^{(2)} &= \begin{bmatrix}
\hat{x}_1 \otimes \hat{x}_1 \\
\hat{x}_2 \otimes \hat{x}_2 \\
\vdots \\
\hat{x}_K \otimes \hat{x}_K
\end{bmatrix} \in \mathbb{R}^{K \times n^2}, \\
R &= \hat{\dot{X}}.
\end{align}

The reduced operators are then obtained by solving the following least-squares problem:
\begin{equation}
\label{eq:lsq-problem}
\min_{O \in \mathbb{R}^{n \times (n + p + n^2)}} \| D O^T - R \|_F^2,
\end{equation}
% where \( O = [\hat{A} \;\; \hat{B} \;\; \hat{F}] \).
where \( O = [\hat{A}^{n \times n} \;\; \hat{B}^{n \times p} \;\; \hat{F}^{n \times n^2}] \).
\subsubsection{Least-Squares Solution via Singular Value Decomposition (SVD)} \label{sec:svd-solve}

The least-squares problem in \eqref{eq:lsq-problem} can be solved robustly using the singular value decomposition (SVD) of the data matrix \( D \). Let:
\[
D = U \Sigma V^T,
\]
where:
\begin{itemize}
    \item \( U \in \mathbb{R}^{K \times r} \), \( V \in \mathbb{R}^{(n + p + n^2) \times r} \) have orthonormal columns,
    \item \( \Sigma \in \mathbb{R}^{r \times r} \) is a diagonal matrix of singular values,
    \item \( r = \operatorname{rank}(D) \).
\end{itemize}

Then the minimum-norm solution is given by:
\begin{equation}
O^T = V \Sigma^{-1} U^T R.
\end{equation}

This corresponds to applying the Moore–Penrose pseudoinverse \( D^\dagger \), i.e., \( O^T = D^\dagger R \).

%\paragraph{Tikhonov Regularization.}
\subsubsection{Regularized Least-Squares: Tikhonov Approach}
In cases where \( D \) is ill-conditioned or nearly rank-deficient, regularization improves stability. The Tikhonov-regularized problem is:
\begin{equation}
\min_{O^T} \| D O^T - R \|_F^2 + \lambda^2 \| O^T \|_F^2,
\end{equation}
where \( \lambda > 0 \) is the regularization parameter.

The solution using SVD is:
\begin{equation}
O^T = \sum_{i=1}^{r} \frac{\sigma_i}{\sigma_i^2 + \lambda^2} (u_i^T R) v_i,
\end{equation}
where \( \sigma_i \) are the singular values, and \( u_i \), \( v_i \) are the corresponding left and right singular vectors.

This regularized solution damps the influence of small singular directions, improving robustness to noise or overfitting.

%%%%%%%%%%%%%%%%%%%%%%%%%%%%%%%%%%%%%%%%%%%%%%%%%
\subsection{Tensor-Train Decomposition and Cross Interpolation}
Tensor-Train (TT) decomposition is a low-rank tensor approximation technique that represents high-dimensional tensors as a sequence of lower-order cores connected via matrix products. It provides an efficient computational framework for compressing and manipulating large-scale data, particularly in problems where the underlying tensor exhibits low-rank structure \cite{oseledets2010tt}.
%\todo[inline]{sections' titles need to be improved !!!}
\subsubsection{Tensor-Train Definition}
\begin{figure}[htbp]
    \begin{center}
    {\large
        \begin{tikzpicture}[scale=2]
            \draw (0,0,0) -- (1,0,0) -- (1,1,0) -- (0,1,0) -- cycle;
            \draw (1,0,0) -- (1,1,0) -- (1,1,-1) -- (1,0,-1) -- cycle;
            \draw (0,1,0) -- (1,1,0) -- (1,1,-1) -- (0,1,-1) -- cycle;
            \node[] at (1.9,0.75,0) {\Large =};
            \draw (2.5,0,0) -- (2.9,0,0) -- (2.9,1,0) -- (2.5,1,0) -- cycle;
            \draw (3.4,1,0) -- (3.8,1,0) -- (3.8,1.3,0) -- (3.4,1.3,0) -- cycle;
            \draw (3.8,1,0) -- (3.8,1.3,0) -- (3.8,1.3,-1) -- (3.8,1,-1) -- cycle;
            \draw (3.4,1.3,0) -- (3.8,1.3,0) -- (3.8,1.3,-1) -- (3.4,1.3,-1) -- cycle;
            \draw (4.6,1,0) -- (5.6,1,0) -- (5.6,1.4,0) -- (4.6,1.4,0) -- cycle;
            \node [] at (-0.2,0.5) {$N_x$};
            \node [] at (0.5,-0.2) {$N_z$};
            \node [] at (1.39,0.1) {$N_y$};
            \node [] at (0.5,0.5) {$\mathcal{X}$};
            \node [] at (2.725,1.15) {$\mathcal{G}_1$};
            \node [] at (2.3,0.5) {$N_x$};
            \node [] at (2.725,-0.15) {$r_1$};
            \node [] at (4,1.85) {$\mathcal{G}_2$};
            \node [] at (3.45,1.55) {$N_y$};
            \node [] at (3.27,1.15) {$r_1$};
            \node [] at (3.65,0.85) {$r_2$};
            \node [] at (5.1,1.55) {$\mathcal{G}_3$};
            \node [] at (5.1,0.83) {$N_z$};
            \node [] at (4.5,1.16) {$r_2$};
            \node[] at (5.9,1.16) {$+\;\;\varepsilon$};
        \end{tikzpicture}}
    \end{center}
    \caption{TT decomposition of a 3-dimensional tensor $\mathcal{X}$ with a truncation error of $\varepsilon$.}
    \label{fig:TT-compression}
\end{figure}
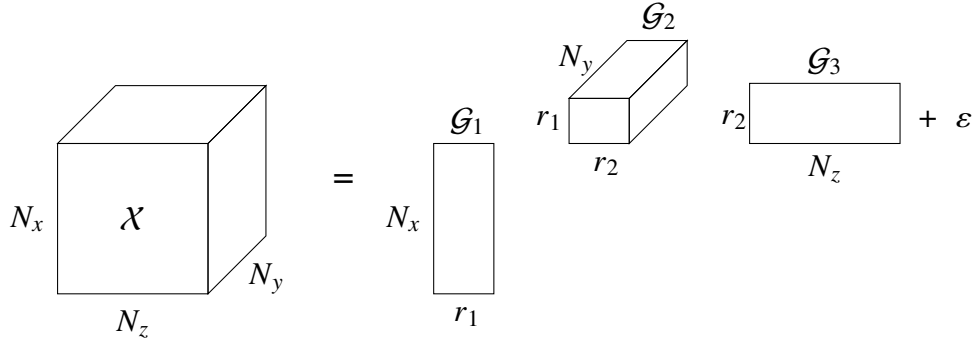

A tensor-train (TT) is an approximate decomposition of a multidimensional tensor \cite{oseledets2011tensor}. For example, a 3-dimensional tensor $\mathcal{X}{(i,j,k)}$ for $1\le i\le N_x$, $1\le j\le N_y$ and $1\le k\le N_z$ can be approximated by the following decomposition
\begin{equation}
    \mathcal{X}(i,j,k)\approx\mathcal{X}_{TT}(i,j,k)=G_1(i)G_2(j)G_3(k)+\varepsilon(i,j,k),
\end{equation}
which is also shown in \Cref{fig:TT-compression}. Here, TT cores are given as $G_1(i)\in\mathbb{R}^{1\times r_1}$, $G_2(j)\in\mathbb{R}^{r_1\times r_2}$ and $G_3(k)\in\mathbb{R}^{r_2\times 1}$, and $(r_1,r_2)$ are known as the TT ranks. Alternatively, the TT decomposition of $X$ can be denoted by $X\approx X_{TT}=G_1\circ G_2\circ G_3$. Note that the TT format of a tensor becomes very efficient for high-dimensional tensors, especially when the TT ranks are significantly smaller than the number of discrete points in each dimension. 

% \subsubsection{TT Arithmetics}
% \todo[inline]{Discuss some TT arithmetics being used in the paper, such as permutation, concatenation, convert from TT-Tensor to TT-matrix, ...}
\subsubsection{TT cross interpolation}
TT cross interpolation is an efficient algorithm for constructing a tensor-train (TT) representation of a high-dimensional tensor without requiring full access to all of its entries. Unlike the TT-SVD algorithm, which operates on the entire tensor and relies on explicit unfolding and singular value decomposition along each mode, TT cross builds the TT cores adaptively by selecting and querying a small subset of tensor entries.

The technique builds upon ideas from skeleton (or CUR) decomposition~\citep{mahoney2009cur} and the Maximum Volume Principle~\citep{goreinov2010find}, which guide the selection of representative tensor slices or fibers. These principles are embedded in heuristic algorithms for TT construction, including Alternating Minimal Energy (AMEn)~\citep{dolgov2014alternating} and Density Matrix Renormalization Group (DMRG)~\citep{savostyanov2011fast} methods.

In the context of this work, TT cross is used to compress the snapshot tensor directly from sampled entries, allowing us to avoid explicit construction of the full snapshot matrix. This is especially valuable for large-scale or unstructured problems where memory constraints or data access patterns preclude dense tensor storage.

\subsection{Tensor Formulation of Full-Order OpInf} \label{sec:full-opinf-tensor}

In the standard OpInf formulation (Section~\ref{sec:opinf-standard}), the snapshot data is stored in a matrix \( X \in \mathbb{R}^{K \times N} \), where each row corresponds to a state vector at a given time. To enable tensor-based compression and operator learning, we reinterpret this data as a structured 5D tensor
\begin{equation}
    x = x_{i,j,k,q,n} \in \mathbb{R}^{N_x \times N_y \times N_z \times N_q \times K},
\end{equation}
where \( (N_x, N_y, N_z) \) are the spatial grid sizes, \( N_q \) is the number of conserved variables (e.g., \( N_q = 3 \) for density, momentum, and energy), and \( K \) is the number of time steps. The total state dimension satisfies \( N = N_x N_y N_z N_q \), matching the structure assumed in Section~\ref{sec:opinf-standard}.

This multidimensional representation allows us to exploit the inherent spatial and physical structure of the data and prepares the model for TT representation later. While the full-order tensor formulation described in this section does not perform model reduction, it serves as a critical foundation for the TT-based OpInf method, which can efficiently handle high-dimensional systems without flattening the data arrays.

Using the 5D snapshot tensor $X$, we will next assemble the data tensor $D$ and time derivative tensor $R$, and set up the least-squares problem in the tensor format.

First, we reshape $x$ into a 9D tensor where the middle four dimensions are singletons to enable tensor products
\begin{equation}
    \begin{aligned}
       X_{i,j,k,q,i',j',k',q',n}&=\text{reshape}(x,[N_x,N_y,N_z,N_{q},1,1,1,1,K])\in\mathbb{R}^{N_x\times N_y\times N_z\times N_{q}\times 1\times 1\times 1\times 1\times K}.
    \end{aligned}
\end{equation}
Second, we construct the quadratic term $X^2$ as the tensor product 
\begin{equation}
\label{eq:quadratic_full}
    X^2=X^2_{i,j,k,q,i',j',k',q',n}=x_{i,j,k,q,n}\,x_{,i',j',k',q',n}\in\mathbb{R}^{N_x\times N_y\times N_z\times N_{q}\times N_x\times N_y\times N_z\times N_{q}\times K}
\end{equation}
where the repeated index $n$ does not imply summation. With this understanding, the data tensor
$$
D_{i,j,k,q,i',j',k',q',n}\in\mathbb{R}^{N_x\times N_y\times N_z\times N_{q}\times (N_x+1)\times (N_y+1)\times (N_z+1)\times (N_{q}+1)\times K}
$$
is defined as
\begin{equation}
    D_{i,j,k,q,i',j',k',q',n}=
    \begin{cases}
        X_{i,j,k,q,i',j',k',q',n} & \text{if } i'=j'=k'=q'=1,\\
        X^2_{i,j,k,q,i'-1,j'-1,k'-1,q'-1,n} & \text{if }i',j',k',q'>1,\\
        0 & \text{otherwise}
    \end{cases}
\end{equation}
Note that the source term $U$ is omitted for simplicity. After assembling the data tensor, we calculate the time derivative tensor 
$$
R_{i,j,k,q,n}\in\mathbb{R}^{N_x\times N_y\times N_z\times N_{q}\times K}
$$
by applying a suitable finite difference formula with desired order of accuracy on the last dimension of $x_{i,j,k,q,n}$.

In the full-order tensorized format, the least-squares problem analogous to \Cref{eq:lsq-problem} is formed as
\begin{equation}\label{eq:ls-full}
    \min_{O}\|{O}_{i,j,k,q,i',j',k',q',i'',j'',k'',q''}D_{i',j',k',q',i'',j'',k'',q'',n}-R_{i,j,k,q,n}\|^2_F
\end{equation}
Notice that the solution tensor $O$ is 12 dimensional tensor such that
$$
O\in\mathbb{R}^{N_x\times N_y\times N_z\times N_{q}\times N_x\times N_y\times N_z\times N_{q}\times (N_x+1)\times (N_y+1)\times (N_z+1)\times (N_{q}+1)}
$$
After \Cref{eq:ls-full} is solved, the full-order operators are extracted from $O$ as shown below
\begin{equation}
    \begin{aligned}
        A_{i,j,k,q,i',j',k',q'}&=O_{i,j,k,q,i',j',k',q,1,1,1,1},\\
        F_{i,j,k,q,i',j',k',q',i''-1,j''-1,k''-1,q''-1}&=O_{i,j,k,q,i',j',k',q',i'',j'',k'',q''}\quad\text{for }1< i'',j'',k'',q''.
    \end{aligned}
\end{equation}
Finally, the prediction is performed by solving the full-order ODE
\begin{equation}\label{eq:full-order-ODE-tensorized}
    \frac{d}{dt}x_{i,j,k,q} = A_{i,j,k,q,i',j',k',q'}x_{i',j',k',q'} + F_{i,j,k,q,i',j',k',q',i'',j'',k'',q''}x_{i',j',k',q'}x_{i'',j'',k'',q''}.
\end{equation}
Notice that \Cref{eq:full-order-ODE-tensorized} is the tensorized version of \Cref{eq:full-ode}, which implies that the full-order OpInf scheme approximates the discretization operators used in the simulation that generated the snapshot data. As mentioned previously, this approach is prohibitively expensive as the number of grids increases in each dimension. However, this is a favorable format for translating this scheme into the full-order TT OpInf method.
% as it is well-known that the TT format is very effective to mitigate the curse of dimensionality \cite{alexandrov2023challenging}. 

Note also that this version of the full-order OpInf method is only applicable when the training dataset is available in the structured grid format. For unstructured grids, the natural tensor structure is lost since spatial coordinates does not follow a regular pattern. The Quantized Tensor-Train (QTT) approach addresses this limitation by using prime factorization to create an artificial tensor structure suitable for TT operations. For example, suppose that the available training dataset is obtained from an unstructured grid simulation with $N$ number of cells, $N_q$ number of conserved variables, and $K$ number of time steps, such that $x\in\mathbb{R}^{N\times N_q\times K}$. The tensorization is then accomplished by considering the prime factorization of $N$ such that $N=N_1N_2\cdots N_m$, where $m$ is the number of the prime factors of $N$, and then, reshaping the snapshot tensor $x$ as
\begin{equation}
    \tilde{x} = \text{reshape}(x,[N_1,N_2,\cdots,N_m,N_q,K])\in\mathbb{R}^{N_1\times N_2\times \cdots\times N_m\times N_q\times K}
\end{equation}
If $m=3$, then the tensorization is the same as the structured grid version described above. For $m\ne3$, a similar tensorization procedure is followed, where the dimensions corresponding to the prime factors are treated like the spatial dimensions of the structured case. This approach forms the basis for the full-order Quantized Tensor-Train (QTT) OpInf method.

\subsection{Full-Order TT and QTT Operator Inference}

Building on the tensor formulation in Section~\Cref{sec:full-opinf-tensor}, we now introduce a compressed operator inference approach using the Tensor-Train (TT) and Quantized Tensor-Train (QTT) formats. In the TT-OpInf framework, the snapshot data, derivative tensors, and learned operator are all represented in TT or QTT format. The least-squares problem is solved directly in this compressed form, avoiding the need to construct large flattened matrices. The QTT variant further extends this approach to unstructured or vectorized data by reshaping the state dimension into a higher-order tensor.

In the following subsections, we describe how the snapshot data is compressed, how the data and derivative tensors are constructed in TT/QTT form, and how the least-squares operator learning problem is solved directly in the tensor-train format.

\subsubsection{Snapshot Tensor and TT/QTT Compression}
As described in Section~\ref{sec:full-opinf-tensor}, the snapshot data is arranged either as a structured 5D tensor for regular grids or as a high-dimensional reshaped tensor for unstructured grids. In both cases, the snapshot tensor is compressed into its TT/QTT format, \(X_{TT}\), to enable efficient storage and computation.

Two compression methods are used: TT-SVD and TT cross interpolation. TT-SVD applies sequential singular value decompositions to compute the TT cores and requires the full snapshot tensor to be stored in memory. In contrast, TT cross interpolation constructs the TT representation adaptively by sampling only selected tensor entries, making it well-suited for large-scale or out-of-core datasets.The TT format of the snapshot tensor forms the basis for all subsequent operations in full-order TT and QTT OpInf.

\subsubsection{TT/QTT-Based Construction of Data and Derivative Tensors}

Once \(X_{TT}\) is available, we construct the data tensor \( D_{TT} \) and the time derivative tensor \( R_{TT} \) required for operator learning.

The TT format of the nonlinear term \( X^2_{TT} \in \mathbb{R}^{n_1 \times \cdots \times n_d \times n_1 \times \cdots \times n_d \times K} \) is computed using TT-cross procedure to approximate the quadratic interaction
\[
X^2(i_1, \ldots, i_d, i'_1, \ldots, i'_d, n) = x(i_1, \ldots, i_d, n) \cdot x(i'_1, \ldots, i'_d, n),
\]
as described earlier in Eq.~\eqref{eq:quadratic_full}. For QTT data, the indices \( i_1, \ldots, i_d \) represent the reshaped spatial modes from the prime factorization.

The data tensor \( D_{TT} \) is then assembled by concatenating the TT representations of the linear and quadratic terms—\( X_{TT} \) and \( X^2_{TT} \)—along a new dimension that mimics column stacking in the matrix formulation. This operation mirrors the structure of the design matrix \( D \) in Section~\ref{sec:opinf-formulation}, and preserves its compositional meaning in compressed form.
% \todo[inline]{How about explaining about the TT-concatenation procedure and and example in TT-Arithmetics. Maybe move that section to an appendix.}

The time derivative tensor \( R_{TT} \in \mathbb{R}^{n_1 \times \cdots \times n_d \times K} \) is computed by applying a finite-difference stencil along the last TT core, which represents the snapshot (temporal) dimension. This operation can be performed directly on the TT cores, preserving both accuracy and rank efficiency. Further implementation details are provided in~\cite{truong2024tensor,manzini2023tensor}.

\subsubsection{TT Least-Squares Operator Learning}

Once the TT format of data tensor \( D_{TT} \) and the time derivative tensor \( R_{TT} \) are constructed, we solve the regularized least-squares problem
\begin{equation}\label{eq:ls-tt}
\min_{O_{TT}} \| D_{TT} O_{TT}^T - R_{TT} \|_F^2 + \gamma \| O_{TT}^T \|_F^2,
\end{equation}
where \( \gamma > 0 \) is a regularization parameter.

This formulation is the TT analog of the matrix problem solved via SVD in Section~\ref{sec:svd-solve}. However, instead of flattening the tensors, all computations are performed directly in the TT format, leveraging low-rank structure throughout.

The solution proceeds in two main steps:
\begin{enumerate}
    \item Compute the regularized TT pseudoinverse of \( D_{TT} \) via core-wise orthogonalization, SVD of the interface core, and TT-permutation.
    \item Apply the TT pseudoinverse to \( R_{TT} \) to construct the TT representation of the learned operator \( O_{TT} \).
\end{enumerate}

The detailed algorithms for TT pseudoinverse and TT least-squares construction are provided in Algorithms~\ref{alg:tt-pinv} and~\ref{alg:tt-ls-solver}, respectively.
\vspace{1ex}
\paragraph{TT Pseudoinverse.}
The first step is to compute the TT representation of the pseudoinverse of the matricized tensor \( D \), where the unfolding is performed at a designated split index. In our implementation, we reshape \( D \) such that the first eight modes correspond to the reduced output dimensions and the remaining modes correspond to the input dimensions. This choice reflects the tensor layout of the operator \( O \in \mathbb{R}^{n_1 \times \cdots \times n_8 \times n_9 \times \cdots \times n_{12}} \), where the left modes represent the output space and the right modes represent the product of reduced coordinates from the snapshot and nonlinear terms.

The TT pseudoinverse is computed using Algorithm~\ref{alg:tt-pinv}, which orthogonalizes the TT cores of \( D \) from both the left and right, computes an SVD of the core interface, applies Tikhonov regularization to the singular values, and then constructs a TT representation of the pseudoinverse \( D^\dagger \) via contraction and transposition of the orthogonal components.
\begin{algorithm}[H]
\DontPrintSemicolon
\KwIn{TT tensor \(\mathcal{X}_{TT}\), split index \(\texttt{coreidx}\), tolerance \(\texttt{tt\_tol}\), regularization \(\gamma\) (optional)}
\KwOut{TT tensor \(\mathcal{X}^\dagger_{TT}\), pseudo-inverse of the matricization of \(\mathcal{X}_{TT}\)}

% \textbf{Step 0:} Set \(\gamma \gets 0\) if not provided\;

\textbf{Step 1:} Extract TT cores \(G \);
% \gets \texttt{core2cell}(\mathcal{X}_{TT})\); store mode sizes \(n\) and TT ranks \(r\)\;

\textbf{Step 2:} Left-orthogonalize the first \(\texttt{coreidx}\) cores:\;
\Indp
\([G_l, R_l] \gets \texttt{tt\_left\_orthogonolize}(G_{1:\texttt{coreidx}})\)\;
\Indm

\textbf{Step 3:} Right-orthogonalize the remaining cores:\;
\Indp
\([G_r, R_r] \gets \texttt{tt\_right\_orthogonolize}(G_{\texttt{coreidx}+1:end})\)\;
\Indm

\textbf{Step 4:} Compute compact SVD: \([U, S, V] \gets \texttt{svd}(R_l R_r)\)\;

\textbf{Step 5:} Apply regularization to singular values:\;
\Indp
\If{\(\gamma > 0\)}{
  \(s \gets \text{diag}(S)\)\;
  \(s \gets \text{diag}(s ./ (s^2 + \gamma))\)\;
}
\Else{
  \(s \gets S^\dagger\) \tcp*{Pseudoinverse of \(S\)}
}
\Indm

\textbf{Step 6:} Update last core of \(G_l\): \(G_l[\text{end}] \gets \texttt{tensorprod}(G_l[\text{end}], U, 3, 1)\)\;

\textbf{Step 7:} Update first core of \(G_r\): \(G_r[1] \gets \texttt{tensorprod}(V, G_r[1], 1, 1)\)\;

\textbf{Step 8:} Assemble the TT tensor \(\mathcal{U}_{TT}\) from the left-side cores \(G_l\), and transpose it by exchanging input and output dimensions in its mode ordering\;

\textbf{Step 9:} Assemble the TT tensor \(\mathcal{V}_{TT}\) from the right-side cores \(G_r\), and transpose it similarly to obtain \(\mathcal{V}_{TT}^T\)\;

\textbf{Step 10:} Scale the final core of \(\mathcal{V}_{TT}^T\) by the regularized singular values \(s\)\;

Combine the transposed tensors \(\mathcal{V}_{TT}^T\) and \(\mathcal{U}_{TT}^T\) to form the TT representation of the pseudo-inverse: \(\mathcal{X}^\dagger_{TT}\)\;

\caption{Pseudo-Inverse of a TT Tensor via Low-Rank SVD and Regularization - \texttt{tt\_pinv}}
\label{alg:tt-pinv}
\end{algorithm}
\vspace{1ex}
\paragraph{Operator Construction.}
After obtaining the TT pseudoinverse \( D^\dagger_{TT} \), the operator tensor \( O_{TT} \) is constructed by contracting \( D^\dagger_{TT} \) with the right-hand side tensor \( R_{TT} \), as outlined in Algorithm~\ref{alg:tt-ls-solver}.

Since, \( D^\dagger_{TT} \) is a $9$-dimensional (4 spatial + 4 operator + 1 time) TT , and the right-hand side \( R_{TT} \) is a $5$-dimensional (4 spatial + 1 time) TT. The contraction eliminates the time dimension and produces the 12D operator tensor matching the dimensional structure of the full operator defined in Section~\ref{sec:full-opinf-tensor}.

The algorithm is written with these prescribed TT dimensions in mind to improve clarity and transparency. While it is tailored to the specific tensor layout used in the full-order TT-OpInf framework, the same approach can be readily generalized to the QTT case. The contraction and assembly steps remain structurally identical, with the only change being the number and interpretation of the TT modes.

\begin{algorithm}[H]
\DontPrintSemicolon
\KwIn{\( D_{TT} \) (data), \( R_{TT} \) (right-hand side), regularization parameter \(\gamma\), TT tolerance \(\varepsilon_{\text{TT}}\)}
\KwOut{ \( O_{TT} \), solution to the regularized least-squares problem}

\textbf{Step 1:} Compute the regularized TT pseudoinverse of \( D_{TT} \):\;
\Indp
\(\tilde{D}_{TT} \gets \texttt{tt\_pinv}(D_{TT},\; \texttt{split\_index}=8,\; \varepsilon_{\text{TT}},\; \gamma)\)\;
\Indm

\textbf{Step 2:} Extract TT cores of \(\tilde{D}_{TT}\) and \( R_{TT} \):\;
\Indp
\([G_D^{(1)}, \dots, G_D^{(9)}] \gets \texttt{cores of } \tilde{D}_{TT}\)\;
\([G_R^{(1)}, \dots, G_R^{(5)}] \gets \texttt{cores of } R_{TT}\)\;
\Indm

\textbf{Step 3:} Initialize TT cores of solution \( O \): \( G_O^{(1)} \gets G_R^{(1)},\; \dots,\; G_O^{(4)} \gets G_R^{(4)} \)\;

\textbf{Step 4:} Contract final core of \( R \) with first core(s) of \( \tilde{D} \):\;
\Indp
\If{\( G_D^{(1)} \) is 2D}{
    \( T \gets \texttt{squeeze}(G_R^{(5)}) \cdot {G_D^{(1)}}^T \)
}
\Else{
    \( T \gets \texttt{squeeze}(G_R^{(5)}) \cdot \texttt{squeeze}(G_D^{(1)}) \)
}
\Indm

\textbf{Step 5:} Set core 5 of \( O \): \( G_O^{(5)} \gets \texttt{tensorprod}(T, G_D^{(2)},\; \text{contract mode } 2 \to 1) \)\;

\textbf{Step 6:} Set cores 6–12 of \( O \): \( G_O^{(6)} \gets G_D^{(3)},\; \dots,\; G_O^{(12)} \gets G_D^{(9)} \)\;

\textbf{Step 7:} Assemble TT tensor \( O \) from its cores:\;
\Indp
\( O \gets \texttt{TT tensor from } [G_O^{(1)}, \dots, G_O^{(12)}] \)\;
\Indm

\caption{TT Least-Squares Solver for Operator Inference}

\label{alg:tt-ls-solver}
\end{algorithm}

% Rather than performing a full contraction across all TT cores, our method builds \( O_{TT} \) by selectively combining the TT cores of \( R_{TT} \) with the first few cores of \( D^\dagger_{TT} \), followed by TT contractions and concatenation with the remaining TT cores of \( D^\dagger_{TT} \). 

\paragraph{Connection to Standard SVD-Based Solution.}
This method serves as the TT-based analog of the SVD approach described in Section~\ref{sec:svd-solve}. Instead of relying on dense matrix operations, the TT formulation performs all computations through core-wise manipulations—such as orthogonalization, low-rank decomposition, and contraction—entirely within the TT format. Tikhonov regularization is supported in both formulations.

By working directly in the compressed TT representation, this approach enables full-order operator inference for high-dimensional systems that would otherwise be intractable using matrix-based methods.
% \vspace{1ex}
\subsubsection{Operator Extraction from TT Representation}

After solving the TT-based least-squares problem, the full operator tensor \( O_{TT} \in \mathbb{R}^{n_1 \times \cdots \times n_{12}} \) contains contributions from both the linear and quadratic terms, encoded along the final four TT cores.

To extract the linear operator \( A_{TT} \), we slice the TT cores corresponding to the input dimensions to isolate the entries associated with the linear term. Specifically, we retain the first slice of the mode along each of the final four TT cores (modes 9–12), corresponding to index 1 in each of those dimensions. The resulting TT tensor contains 8 output modes and 4 singleton input modes, which are removed through reshaping:
\[
A_{TT} \in \mathbb{R}^{N_x \times N_y \times N_z \times N_q \times N_x \times N_y \times N_z \times N_q}.
\]

Similarly, to extract the quadratic operator \( F_{TT} \), we slice the same TT cores to remove the first index and retain all remaining slices, which represent the quadratic interaction modes. The resulting TT has full input dimensions:
\[
F_{TT} \in \mathbb{R}^{N_x \times N_y \times N_z \times N_q \times N_x \times N_y \times N_z \times N_q \times N_x \times N_y \times N_z \times N_q}.
\]

\subsubsection{Prediction and Reconstruction in TT Format}

Once the operator tensors \( A_{TT} \) and \( F_{TT} \) have been learned, the full-order solution in TT format is advanced in time by solving the nonlinear ODE system directly in TT format. The TT-encoded initial condition \( x(t=0) \) is first rounded to a prescribed truncation tolerance \( \varepsilon_{TT} \) to ensure rank efficiency.

The time integration in TT format is performed using an adaptive Runge–Kutta method RK45. At each timestep, the right-hand side function is evaluated via TT-contractions involving the learned operators \( A_{TT} \) and \( F_{TT} \), and optionally \( B_{TT} \) if source terms are included. After solving the system over the interval \( [0, T_{\text{final}}] \), the final state \( x(T_{\text{final}}) \) remains in TT format and can be post-processed or reshaped as needed. 
%%%%%%%%%%%%%%%%%%%%%%%%

\subsection{Reduced-order TT-OpInf}
In this section, we present the reduced-order version of the TT-OpInf method. The goal of this approach is to eliminate the singular value decomposition (SVD) step used in the standard OpInf method \cite{peherstorfer2016opinf,kramer2024review} to compute the POD basis matrix $V$. Instead, the snapshot tensor $X \in\mathbb{R}^{N_x \times N_y \times N_z \times N_q\times K}$ is directly converted to the TT format 
\begin{equation}
    X_{TT}=G_1 G_2 G_3 G_4 G_5
\end{equation}
either by TT-SVD \cite{oseledets2011tensor} or TT cross interpolation \cite{oseledets2010tt} algorithms, where $G_i$ are the TT cores. Interestingly, with $G_i$ for $i=1,\dots,4$ being orthogonalized, the reduced snapshot matrix $\hat{X}$ in the standard OpInf method \cite{peherstorfer2016opinf} is equivalent to $G_5$, such that $\hat{X}=G_5$. Since $\hat{X}$ is available as the second core of $X_{TT}$ by default after the TT conversion is completed, it can be used to learn the reduced operators $\hat{A},\hat{F},\hat{B}$ exactly as described in \Cref{sec:opinf-standard}. After the solving \Cref{eq:reduced-ode}, the reduced-order predicted solution $\hat{x}$ can be transformed back to the full-order space by simply forming a new TT
\begin{equation}
    x_{TT}=G_1 G_2 G_3 G_4\hat{x}
\end{equation}
The key advantage of this approach is that it replaces the traditional POD projection via SVD with a direct TT conversion of the snapshot tensor. In particular, when TT cross interpolation is used, the snapshot data can be compressed without ever forming the full tensor in memory. Moreover, the resulting TT representation not only identifies a reduced coordinate system through the final TT core (e.g., \( \hat{X} = G_5 \)), but also compresses the full snapshot tensor along all dimensions. This compression leads to significant memory and computational savings, making the reduced-order TT-OpInf method practical for problems involving very large spatial grids or long time horizons. As a result, TT-ROM enables operator learning and prediction in settings that would otherwise be prohibitively expensive using standard matrix-based approaches.

%%%%%%%%%%%%%%%%%%%%%%%%%%%%%%%%%%%%%%%%%%%%%%%%%%%%%%%%%%%%
\section{Results}

In this section, we evaluate the performance of various operator inference methods—including standard ROM-OpInf, full-order tensorized (FT) OpInf, tensor-train (TT) and quantized tensor-train (QTT) formulations, as well as the proposed reduced-order TT-ROM approach—on three benchmark problems. These include the linear 2D heat equation, the nonlinear 2D viscous Burgers' equation, and compressible laminar flow over a cylinder on an unstructured mesh.

For each test case, we report the computational cost of three key phases: dimensionality reduction or TT compression (\( t_{\text{POD}} \)), operator learning (\( t_{\text{learn}} \)), and prediction (\( t_{\text{predict}} \)). Accuracy is assessed in terms of relative error with respect to high-fidelity simulation data. All TT-based methods use a specified truncation tolerance \( \varepsilon_{\text{TT}} \), and regularization parameters are set appropriately for each method.

These benchmarks highlight the scalability, accuracy, and computational trade-offs of TT-based methods relative to standard and full-tensor OpInf approaches.
\subsection{Heat Equation}
In this example, we consider the 2-dimensional heat equation on the computational domain $\Omega$,
\begin{equation}\label{eq:heat2d}
    u_t=\mu\nabla^2u,\quad\Omega\in[0,1]\times[0,1]
\end{equation}
where $\mu=0.01$ the diffusion coefficient. The data is generated by solving \Cref{eq:heat2d} using a $4^{th}$-order conservative central finite difference scheme until the final time $T=1$ on a uniform $20\times20$ grid. The initial condition at $t=0$ is given by 
\begin{equation}
    u_0(x,y)=\cos{2\pi(x+y)}
\end{equation}
and periodic boundary conditions are applied. All OpInf schemes use snapshots generated up to $t = 0.25$ for training, and report predictions evaluated at $t = 0.9$. Note that this example only tests the linear operator of the OpInf framework due to the lack of the nonlinear convection term. 

In this problem, both the ROM and TT-ROM methods use a regularization factor of $10^{-6}$ for the linear operator $A$ and the nonlinear operator $F$. In contrast, the full-order TT and QTT OpInf schemes do not require any regularization. Additionally, a very small truncation tolerance of $\varepsilon_{\text{TT}} = 10^{-12}$ is used for all TT-based computations. The standard ROM is executed with the same number of POD modes $n$ as obtained from the TT-ROM under this truncation tolerance.

The performance comparison for various OpInf variants is presented in \Cref{tab:heat2d}. The ROM method corresponds to the standard OpInf approach \cite{peherstorfer2016opinf,kramer2024review}. TT-ROM refers to a reduced-order OpInf method in which the snapshot data is first converted to tensor-train (TT) format to define reduced variables, followed by learning and conversion back to the full-order space using the tensor-train method. FT denotes the full-order OpInf method applied directly in full-tensor format, as discussed previously. TT is the tensor-train formulation of FT, and QTT represents the quantized tensor-train version of TT.

\begin{table}[htbp]
\centering
\caption{Performance comparison for the 2D heat problem.}
\begin{tabular}{|l|c|c|c|c|c|}
\hline
\textbf{OpInf} & \boldmath$t_{\text{POD}}$ \textbf{(s)} & \boldmath$t_{\text{learn}}$ \textbf{(s)} & \boldmath$t_{\text{predict}}$ \textbf{(s)} & \textbf{Rel. Error} & \textbf{POD: $n$ or TT-rank: $r$} \\
\hline
ROM   & 3.27e--02 & 2.45e--03 & 4.16e--02 & 1.61e--03  & $n = 1$ \\
TT-ROM& 2.08e--02 & 1.84e--03 & 1.63e--02 & 1.61e--03  & $n = 1$ \\
FT    &       --  & 6.00e--02 & 7.05e--02 & 1.61e--03  & -- \\
TT    & 2.83e--02 & 4.29e--02 & 4.63e--01 & 1.61e--03  & $r = (2, 1, 1)$ \\
QTT   & 4.95e--01 & 4.47e--02 & 1.08e--01 & 1.61e--03  & $r = (2, 2, 4, 8, 2, 1)$ \\
\hline
\end{tabular}
\label{tab:heat2d}
\end{table}

Performance is evaluated based on several metrics: the time required to reduce the snapshot data to a lower-dimensional space ($t_{\text{POD}}$), the time required for learning ($t_{\text{learn}}$), the time required for prediction ($t_{\text{predict}}$), and the maximum relative error (Rel. Error) with respect to the simulation data. The final column reports either the number of POD dimensions $n$ or the tensor-train (TT) ranks $r$, depending on the method.

Note that the TT method operates on a 5-dimensional tensor $X_{i,j,k,q,n}$, where $(i,j,k)$ correspond to the structured spatial grid in the $x$, $y$, and $z$ directions; $q$ is the equation index (e.g., for systems like the Navier–Stokes equations, though $q=1$ for the heat equation); and $n$ indexes time snapshots. For TT ranks, only the internal ranks $r_1$, $r_2$, and $r_3$ are reported, as the boundary ranks $r_0$ and $r_4$ are always equal to 1 by definition and are therefore omitted. Similarly, only the internal ranks of the QTT method are reported.

For the Heat Equation, \Cref{tab:heat2d} shows that all methods achieve the same relative error. The time required to project the snapshot data to a lower-dimensional space is of the same order for the ROM, TT-ROM, and TT method, whereas the QTT method takes nearly an order of magnitude longer. The learning times for the reduced models (ROM and TT-ROM) are comparable, while all full-order models exhibit learning times nearly an order of magnitude higher. During the prediction stage, the TT-based full-order methods are more expensive, primarily due to the repeated rounding operations required within the explicit RK45 scheme to control TT-rank growth.

\subsection{2D Viscous Burgers' Equation}
In this example, we include the nonlinear convection term and consider the 2-dimensional viscous Burgers' equation on the computational domain $\Omega$,
\begin{equation}\label{eq:visc-burgers-2d}
    u_t+uu_x+uu_y=\nabla^2u,\quad\Omega\in[0,4]\times[0,4]
\end{equation}
The data is generated by solving \Cref{eq:visc-burgers-2d} using a $5^{th}$-order conservative upwind finite difference scheme for convective fluxes and a $4^{th}$-order conservative central finite difference scheme for the diffusive fluxes until the final time $T=0.1$ on a uniform $20\times20$ grid. The initial condition at $t=0$ is given by 
\begin{equation}
    u_0(x,y)=\frac{\pi\sin\left(\frac{\pi(x+y)}{2}\right)}{2+\cos\left(\frac{\pi(x+y)}{2}\right)}
\end{equation}
and periodic boundary conditions are applied. All OpInf schemes use snapshots generated up to $t = 0.05$ for training, and report predictions evaluated at $t = 0.09$. 

In this problem, both the ROM and TT-ROM methods employ a regularization factor of $10^{-6}$ for the linear operator $A$ and the nonlinear operator $F$. As before, the full-order TT and QTT OpInf schemes do not require any regularization. A truncation tolerance of $\varepsilon_{\text{TT}} = 10^{-5}$ is used for all TT-based computations. Consistent with the previous case, the standard ROM is executed using the same number of POD modes $n$ as determined by the TT-ROM under this truncation tolerance.

\Cref{tab:viscburgers2d} presents the performance comparison of various OpInf methods for the 2D viscous Burgers' equation. As in the previous case, the QTT method exhibits the largest $t_{\text{POD}}$, while the other methods require similar times for projecting snapshots to the reduced space. The ROM and TT-ROM methods yield comparable $t_{\text{learn}}$ values, as they differ only in the initial data projection step.

In this example, the full-order FT method shows a significantly larger $t_{\text{learn}}$—several orders of magnitude higher—due to the presence of the quadratic term in the least-squares problem and the associated Kronecker product operations. The full-order TT and QTT methods exhibit similar learning times to each other.

In terms of prediction, ROM and TT-ROM perform similarly to the heat equation case. Among the full-order models, TT-OpInf has the highest $t_{\text{predict}}$, likely due to the TT rounding operations required to control rank growth. Despite the QTT method introducing more dimensions through prime factorization, it achieves better performance in both learning and prediction compared to TT-OpInf.
\begin{table}[htbp]
\centering
\caption{Performance comparison for the 2D viscous Burgers' equation.}
\begin{tabular}{|l|c|c|c|c|c|}
\hline
\textbf{OpInf} & \boldmath$t_{\text{POD}}$ \textbf{(s)} & \boldmath$t_{\text{learn}}$ \textbf{(s)} & \boldmath$t_{\text{predict}}$ \textbf{(s)} & \textbf{Rel. Error} & \textbf{POD: $n$ or TT-rank: $r$} \\
\hline
ROM    & 4.31e--02 & 8.25e--03 & 5.25e--02 & 1.09e--03 & $n = 5$ \\
TT-ROM & 3.29e--02 & 9.23e--03 & 3.06e--02 & 1.24e--03 & $n = 5$ \\
FT     &        -- & 3.98e+01  & 9.35e--01 & 1.09e--03 & -- \\
TT     & 7.52e--02 & 3.22e--01 & 1.38e+00 & 1.44e--03 & $r = (18, 1, 1)$ \\
QTT    & 3.34e--01 & 1.43e--01 & 2.65e--01 & 1.25e--03 & $r = (2, 4, 8, 16, 5, 1)$ \\
\hline
\end{tabular}
\label{tab:viscburgers2d}
\end{table}

\subsection{Compressible Laminar Flow over a Cylinder}In this test case, we consider a 2D unsteady laminar flow over a cylinder. The Reynolds number is set to $Re = 75$, which exceeds the critical value for the onset of vortex shedding in the cylinder’s wake. The freestream Mach number is $M = 0.2$, corresponding to the incompressible regime. The problem setup—including the geometry, mesh, and initial and boundary conditions—is implemented exactly as described in \cite{danis2022ddg}. For a complete description of the governing equations and the numerical method, the reader is referred to \cite{danis2022ddg}.

Snapshots are generated using cell-averaged solutions computed from a fifth-order numerical solution obtained with the direct Discontinuous Galerkin method with interface correction (DDGIC) solver introduced in \cite{danis2022ddg}. The mesh used in this problem is an unstructured triangular mesh, consisting of 3,280 cells, and 750 snapshots are collected. The training phase uses data up to $t = 10$ seconds, while predictions are evaluated at $t = 15$ seconds.

In this problem, both the ROM and TT-ROM methods employ a regularization factor of $10^{5}$ for the linear operator $A$ and the nonlinear operator $F$. Unlike the previous cases, the full-order TT and QTT OpInf schemes require a regularization factor of $5 \times 10^9$. A truncation tolerance of $\varepsilon_{\text{TT}} = 10^{-4}$ is used for all TT-based computations. As in the previous case, the standard ROM is executed using the same number of POD modes $n$ as determined by the TT-ROM under this truncation tolerance.

\Cref{tab:cylinderflow} summarizes the performance metrics of the OpInf methods for the 2D cylinder flow case. Due to the unstructured nature of the mesh, only the ROM, TT-ROM, and QTT OpInf methods are evaluated. Unlike the previous test cases, the $t_{\text{POD}}$ for TT-ROM is approximately twice that of the standard ROM. This outcome is expected: TT-ROM employs TT cross-interpolation to compute the reduced basis, whereas standard ROM relies on singular value decomposition (SVD). MATLAB provides a highly optimized implementation of SVD, while TT cross-interpolation is an iterative procedure. The increased $t_{\text{POD}}$ indicates that TT cross-interpolation required more iterations in this case compared to the previous ones. Still, this is not a discouraging outcome. Recall that our TT cross implementation does not require loading the complete snapshot data into memory. Instead, only the necessary portions are accessed, as determined dynamically by the cross interpolation logic. Therefore, the success of our cross interpolation approach in this test case is encouraging for its use in more realistic, large-scale problems. Alternatively, we could have employed the standard TT compression algorithm, which is also based on SVD and is thus quite similar in nature to the POD algorithm used in the standard ROM-OpInf method. In that case, the TT-ROM approach would result in a $t_{\text{POD}}$ very similar to that of the standard ROM-OpInf method. That result is not reported here to avoid duplication.

In addition to these reduced-order models, the full-order QTT OpInf method incurs significantly higher computational cost during the initial conversion of data into QTT format—nearly an order of magnitude greater. As a full-order model, it also requires substantially more time for both learning and prediction, which is also expected.

\begin{table}[htbp]
\centering
\caption{Performance comparison for the 2D cylinder flow case.}
\begin{tabular}{|l|c|c|c|c|}
\hline
\textbf{OpInf} & \boldmath$t_{\text{POD}}$ \textbf{(s)} & \boldmath$t_{\text{learn}}$ \textbf{(s)} & \boldmath$t_{\text{predict}}$ \textbf{(s)} & \textbf{POD: $n$ or TT-rank: $r$} \\
\hline
ROM    & 8.24e--01 & 1.93e--01  & 6.06e--02  & $n = 33$ \\
TT-ROM & 1.84e+00  & 1.91e--01  & 2.18e--02  & $n = 33$ \\
QTT    & 1.93e+01  & 2.37e+00  & 1.62e+00  & $r = (2,\ 4,\ 8,\ 16,\ 80,\ 4)$ \\
\hline
\end{tabular}
\label{tab:cylinderflow}
\end{table}

\begin{figure}[htbp]
    \centering
    \begin{subfigure}[t]{0.48\textwidth}
        \centering
        \includegraphics[width=0.75\textwidth,trim={7cm 7cm 7cm 7cm},clip]{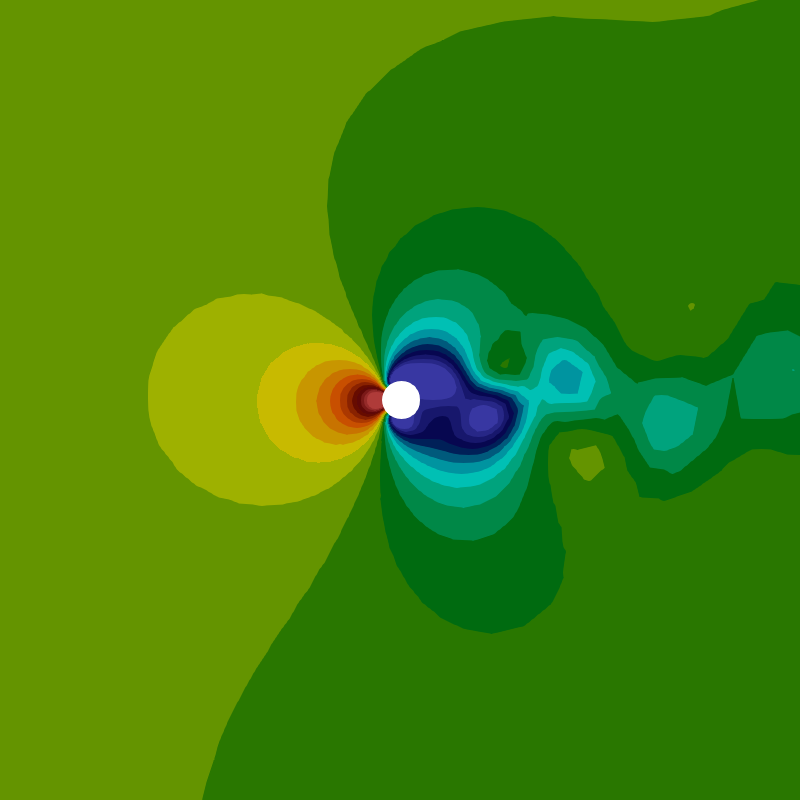}
        \caption{High-fidelity DDGIC simulation}
        \label{fig:sim}
    \end{subfigure}
    \begin{subfigure}[t]{0.48\textwidth}
        \centering
        \includegraphics[width=0.75\textwidth,trim={7cm 7cm 7cm 7cm},clip]{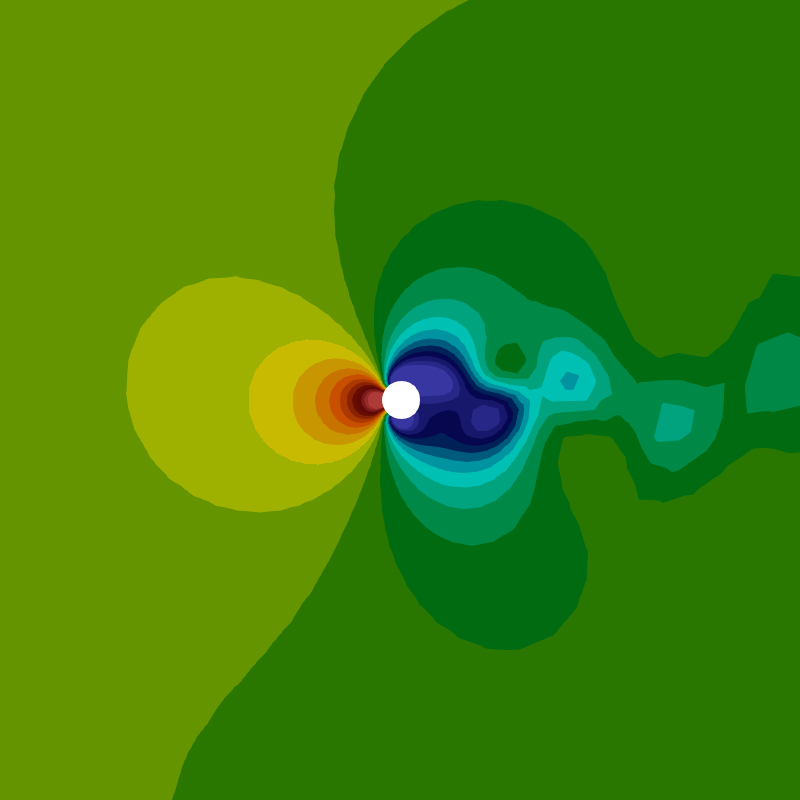}
        \caption{ROM}
        \label{fig:rom}
    \end{subfigure}\vspace{.5cm}
    \begin{subfigure}[t]{0.48\textwidth}
        \centering
        \includegraphics[width=0.75\textwidth,trim={7cm 7cm 7cm 7cm},clip]{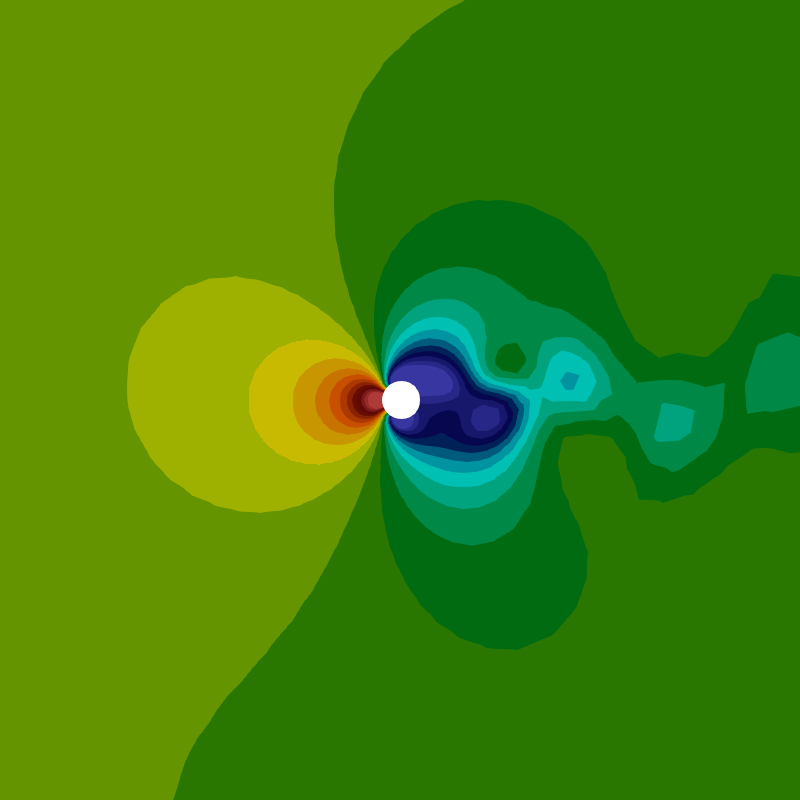}
        \caption{TT-ROM}
        \label{fig:ttrom}
    \end{subfigure}
    \begin{subfigure}[t]{0.48\textwidth}
        \centering
        \includegraphics[width=0.75\textwidth,trim={7cm 7cm 7cm 7cm},clip]{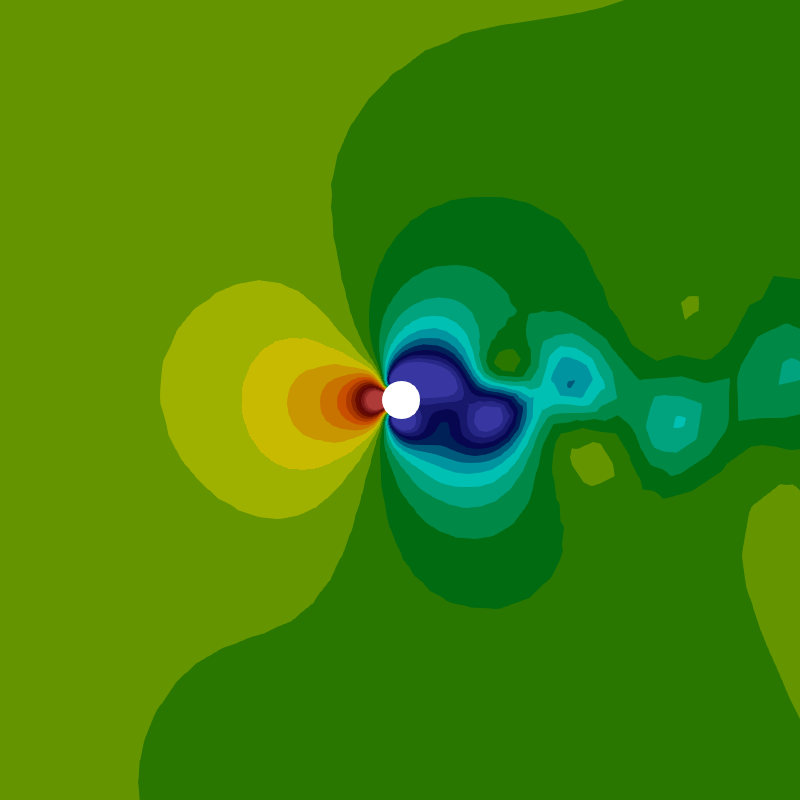}
        \caption{QTT-OpInf}
        \label{fig:qtt}
    \end{subfigure}
    \caption{Comparison of the predicted cell-averaged pressure field at $t = 15$\,s for different OpInf methods. Pressure contours are shown using 25 equispaced levels in the range $[0.70,\ 0.73]$.}
    \label{fig:flow-over-a-cylinder}
\end{figure}

\Cref{fig:flow-over-a-cylinder} shows the pressure field obtained from the cell-averaged DDGIC solution alongside the predictions from the OpInf methods. Compared to the high-fidelity simulation, all OpInf schemes successfully predict the flow beyond the training horizon. Notably, the pressure fields produced by ROM and TT-ROM are nearly identical, while the QTT method exhibits only minor differences.

Across all test problems, the proposed TT and QTT OpInf methods demonstrate comparable accuracy to standard ROM and FT approaches, with significantly improved memory efficiency. The TT-ROM variant achieves competitive training and prediction speed while eliminating the need for full snapshot storage. QTT is particularly effective for unstructured meshes, although it incurs higher compression cost due to mode expansion. These results confirm the scalability and robustness of the TT-OpInf framework, particularly in settings where classical projection-based methods become impractical.

%%%%%%%%%%%%%%%%%%%%%%%%%%%%%%%%
\section{Conclusion}
In this work, we have presented a comprehensive tensor-based formulation of the Operator Inference (OpInf) framework, resulting in a formulation that can handle larger problems with reduced memory overhead, using Tensor-Train (TT) and Quantized Tensor-Train (QTT) formats. The proposed methods extend the applicability of OpInf to extremely high-dimensional dynamical systems, offering new possibilities for operator learning in both structured and unstructured settings.

By leveraging low-rank tensor decompositions, TT and QTT OpInf formulations preserve the multi-dimensional structure of the data and significantly reduce computational and storage costs. The reduced-order TT-OpInf method further eliminates the reliance on SVD-based projection, enabling data-driven modeling even when memory constraints or data scale would make traditional approaches infeasible.

Validation on canonical partial differential equations problems—such as the heat equation and viscous Burgers’ equation—demonstrates the accuracy, efficiency, and flexibility of these methods. While there remain challenges in terms of complexity, implementation robustness, and generalization to more complex geometries or boundary conditions, we believe this work represents a meaningful step toward scalable data-driven modeling for high-dimensional physical systems.

\section*{Acknowledgments}
The authors gratefully acknowledge the support of the Laboratory Directed Research and Development (LDRD) program of Los Alamos National Laboratory under project number 20230067DR.
Los Alamos National Laboratory is operated by Triad National Security, LLC, for the National Nuclear Security Administration of U.S. Department of Energy (Contract No.\ 89233218CNA000001).
\clearpage
\bibliography{references}

\end{document}